\documentclass{amsart}
\usepackage{latexsym, amssymb}
\newtheorem{thm}{Theorem}

\newtheorem{cor}{Corollary}

\newtheorem{lem}{Lemma}
\newtheorem{defn}{Definition}

\begin{document}
\begin{center}
{\bf FINITE TYPE LINK CONCORDANCE INVARIANTS}\\
\vspace{.2in}
{\footnotesize BLAKE MELLOR}\\
{\footnotesize Honors College}\\
{\footnotesize Florida Atlantic University}\\
{\footnotesize 5353 Parkside Drive}\\
{\footnotesize Jupiter, FL  33458}\\
{\footnotesize\it  bmellor@fau.edu}\\
\vspace{1in}
{\footnotesize ABSTRACT}\\
{\ }\\
\parbox{4.5in}{\footnotesize \ \ \ \ \ This paper is a follow-up to \cite{me1}, in which
the author showed that the only real-valued finite type invariants of link
homotopy are the linking numbers of the components.  In this paper, we extend the methods
used to show that the only real-valued finite type invariants of link concordance are, again,
the linking numbers of the components.
\noindent {\it Keywords:}  Finite type invariants; link concordance.}\\
\end{center}
\input{vpsfig.sty}
\tableofcontents
\section{Introduction} \label{S:intro}
This paper extend the arguments in the author's previous work on link homotopy \cite{me1}
to link concordance.  Many of the arguments are nearly identical; the main differences
lie in the base cases to the main theorem, in Section~\ref{S:size}.
We will begin with a brief overview of finite type invariants.
In 1990, V.A. Vassiliev introduced the idea of {\it Vassiliev} or {\it finite type} knot
invariants, by looking at certain groups associated with the cohomology of the space of
knots.  Shortly thereafter, Birman and Lin~\cite{bl} gave a combinatorial description
of finite type invariants.  We will give a summary of this combinatorial theory.
For more details, see Bar-Natan~\cite{bn1}.
\subsection{Singular Knots and Chord Diagrams} \label{SS:chord}
We first note that we can extend any knot invariant to an invariant of
{\it singular} knots, where a singular knot is an immersion of $S^1$ in
3-space which is an embedding except for a finite number of isolated double
points.  Given a knot invariant $v$, we extend it via the relation:
$$\psfig{file=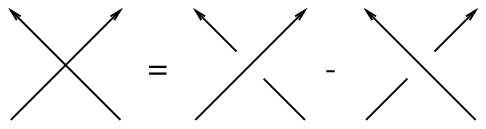}$$
An invariant $v$ of singular knots is then said to be of {\it finite type}, if
there is an integer $d$ such that $v$ is zero on any knot with more than $d$
double points.  $v$ is then said to be of {\it type} $d$.  We denote by $V_d$
the space generated by finite type invariants of type $d$.  We
can completely understand the space of finite type invariants by understanding
all of the vector spaces $V_d/V_{d-1}$.  An element of this vector space is
completely determined by its behavior on knots with exactly $d$ singular
points.  Since such an element is zero on knots with more than
$d$ singular points, any other (non-singular) crossing of the knot can be
changed without affecting the value of the invariant.  This means that
elements of $V_d/V_{d-1}$ can be viewed as functionals on the space of
{\it chord diagrams}:
\begin{defn}
A {\bf chord diagram of degree d} is an oriented circle, together with $d$
chords of the circles, such that all of the $2d$ endpoints of the chords are
distinct.  The circle represents a knot, the endpoints of a chord represent
2 points identified by the immersion of this knot into 3-space.
\end{defn}
Functionals on the space of chord diagrams which are derived from knot
invariants will satisfy certain relations.  This leads us to the definition
of a {\it weight system}:
\begin{defn}
A {\bf weight system of degree d} is a function $W$ on the space of chord
diagrams of degree $d$ (with values in an associative commutative ring
${\bf K}$ with unity) which satisfies 2 relations:
\begin{itemize}
    \item (1-term relation) $$\psfig{file=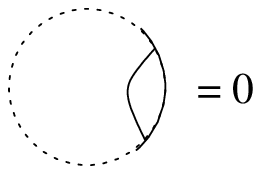}$$
    \item (4-term relation) $$\psfig{file=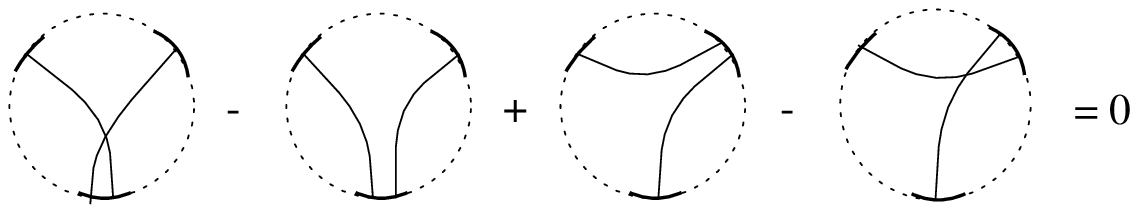}$$
    Outside of the solid arcs on the circle, the diagrams can be anything,
as long as it is the same for all four diagrams.
\end{itemize}
We let $W_d$ denote the space of weight systems of degree d.
\end{defn}
Bar-Natan~\cite{bn1} defines maps $w_d: V_d \rightarrow W_d$ and $v_d: W_d \rightarrow V_d$.
$w_d$ is defined by {\it embedding} a chord
diagram $D$ in ${\bf R}^3$ as a singular knot $K_D$, with the chords corresponding to
singularities of the embedding (so there are $d$ singularities).  Any two
such embeddings will differ by crossing changes, but these changes will not effect the value of
a type $d$ Vassiliev invariant on the singular knot.  Then, for any type $d$ invariant $\gamma$,
we define $w_d(\gamma)(D) = \gamma(K_D)$.  Bar-Natan shows that this is, in fact, a weight
system.  The 1-term relation is satisfied because of the first Reidemeister move, and the 4-term
relation is essentially the result of rotating a third strand a full turn around a double point.
$v_d$ is much more complicated to define, using the Kontsevich integral.  For a full treatment
of the Kontsevich integral, see Bar-Natan~\cite{bn1} and Le and Murakami~\cite{lm}.  Using a Morse
function, any knot (or link or string link) can be decomposed into elementary {\it tangles}:
$$\psfig{file=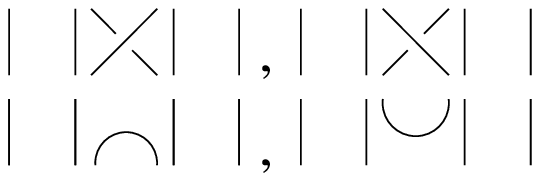}$$
Le and Murakami define a map $Z$ from an elementary tangle with $k$ strands
to the space of chord diagrams on $k$ strands.  This map respects composition of tangles:
if $T_1\cdot T_2$ is the tangle obtained by placing $T_1$ on top of $T_2$, then
$Z(T_1\cdot T_2) = Z(T_1)Z(T_2)$.  Le and Murakami prove that this map gives an isotopy
invariant of knots and links.
Given a degree $d$ weight system $W$, and a knot $K$, we now define $v_d(W)(K) = W(Z(K))$.
Bar Natan shows that $w_d$ and $v_d$ are ``almost'' inverses.  More precisely,
$w_d(v_d(W)) = W$ and $v_d(w_d(\gamma))-\gamma$ is a knot invariant of type $d-1$.  As a result,
(see \cite{bl,bn1,st}) the space $W_d$ of weight
systems of degree $d$ is isomorphic to $V_d/V_{d-1}$.  For convenience,
we will usually take the dual approach, and simply study the space of
chord diagrams of degree $d$ modulo the 1-term and 4-term relations.
The dimensions of these spaces have been computed for $d \leq
12$ (see Bar-Natan~\cite{bn1} and Kneissler~\cite{kn}).  It is useful to combine
all of these spaces into a graded module via direct sum.  We can give this module a Hopf
algebra structure by defining an appropriate product and co-product:
\begin{itemize}
    \item  We define the product $D_1 \cdot D_2$ of two chord diagrams
$D_1$ and $D_2$ as their connect sum.  This is well-defined modulo the 4-term
relation (see \cite{bn1}).
$$\psfig{file=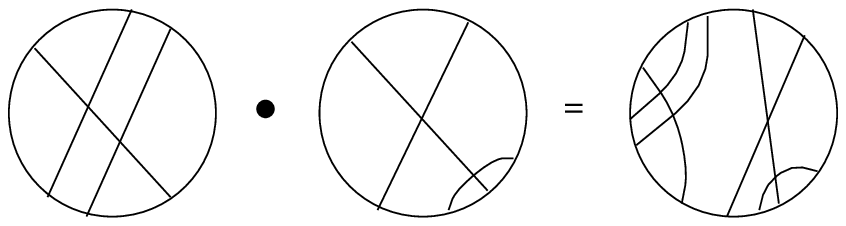}$$
    \item  We define the co-product $\Delta(D)$ of a chord diagram $D$ as
follows:
$$
\Delta(D) = {\sum_J D_J' \otimes D_J''}
$$
where $J$ is a subset of the set of chords of $D$, $D_J'$ is $D$ with all the
chords in $J$ removed, and $D_J''$ is $D$ with all the chords not in J
removed.
$$\psfig{file=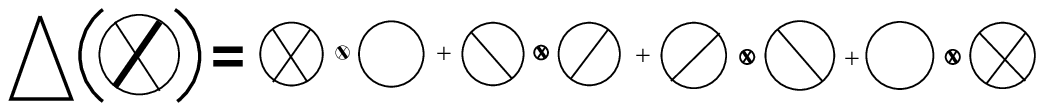}$$
\end{itemize}
It is easy to check the compatibility condition $\Delta(D_1\cdot D_2) = \Delta
(D_1)\cdot\Delta(D_2)$.
\subsection{Unitrivalent Diagrams} \label{SS:ud}
\begin{figure}
    $$\psfig{file=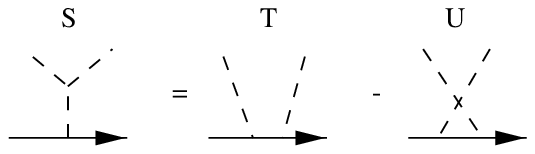}$$
    \caption{STU relation} \label{F:stu}
    \end{figure}

\begin{figure}
    $$\psfig{file=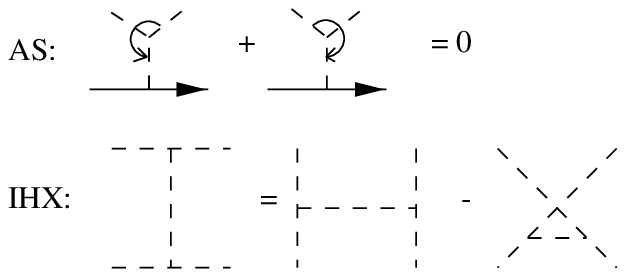}$$
    \caption{AS and IHX relations} \label{F:ihx}
    \end{figure}

It is often useful to consider the Hopf algebra of bounded unitrivalent diagrams, rather
than chord diagrams.  These diagrams, introduced by Bar-Natan~\cite{bn1} (Bar-Natan
calls them {\it Chinese Character Diagrams}), can
be thought of as a shorthand for writing certain linear combinations of chord
diagrams.  We define a {\it bounded unitrivalent graph} to be a unitrivalent graph, with
oriented vertices, together with a bounding circle to which all the univalent vertices are
attached.  We also require that each component of the graph have at least one univalent
vertex (so every component is connected to the boundary circle).  We define the space $A$
of bounded unitrivalent diagrams as the quotient of the space of all bounded unitrivalent
graphs by the $STU$ relation, shown in Figure~\ref{F:stu}.
As consequences of $STU$ relation, the anti-symmetry ($AS$) and $IHX$ relations, see
Figure~\ref{F:ihx}, also hold in $A$.
Bar-Natan shows that $A$ is isomorphic to the algebra of chord diagrams.
We can get an algebra $B$ of {\it unitrivalent diagrams} by simply removing the bounding circle
from the diagrams in $A$, leaving graphs with trivalent and univalent vertices, modulo
the $AS$ and $IHX$ relations.  Bar-Natan shows that the spaces $A$ and $B$ are
isomorphic.  The map $\chi$ from $B$ to $A$ takes a diagram to the linear combination of
all ways of attaching the univalent vertices to a bounding circle, divided by total number of
such ways (T. Le noticed that this factor, missing in \cite{bn1}, is necessary to preserve the
comultiplicative structure of the algebras).  The inverse map $\sigma$
turns a diagram into a linear combination of diagrams by performing sequences of ``basic
operations,'' and then removes the bounding circle.  The two basic operations are:
$$\psfig{file=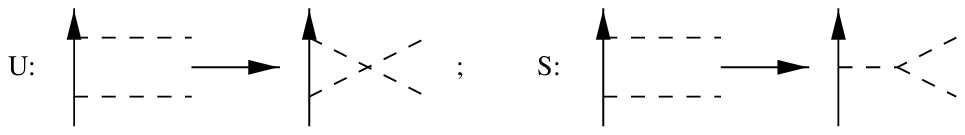}$$
\section{String Links, Links and Concordance} \label{S:string}
\subsection{String Links} \label{SS:string_links}
Bar-Natan~\cite{bn2} extends the theory of finite type invariants to string links.
\begin{defn} (see \cite{hl})
Let D be the unit disk in the plane and let I = [0,1] be the unit interval.  Choose k points
$p_1,..., p_k$ in the interior of D, aligned in order along the the x-axis.  A {\bf string
link} $\sigma$ of k components is a smooth proper imbedding of k disjoint copies of I into
$D \times I$:
$$\sigma:\ \bigsqcup_{i=1}^k{I_i} \rightarrow D \times I$$
such that $\sigma|_{I_i}(0) = p_i \times 0$ and $\sigma|_{I_i}(1) = p_i \times 1$.  The image
of $I_i$ is called the ith string of the string link $\sigma$.
\end{defn}
Essentially, everything works the same way for string links as for knots.
The bounding circle of the bounded unitrivalent diagrams now becomes a set of
bounding line segments, each labeled with a color, to give an algebra $A^{sl}$ (the multiplication
is given by placing one diagram on top of another).  The univalent diagrams are unchanged,
except that each univalent vertex is also labeled with a color to
give the space $B^{sl}$.  The isomorphisms $\chi$ and $\sigma$ between $A$ and $B$ easily
extend to isomorphisms $\chi^{sl}$ and $\sigma^{sl}$ between $A^{sl}$ and $B^{sl}$, just
working with each color separately.  In addition, there are obvious maps $w_d^{sl}$ and
$v_d^{sl}$ analogous to $w_d$ and $v_d$ (we just need to keep track of colors).
\subsection{Links} \label{SS:links}
\begin{figure}
    $$\psfig{file=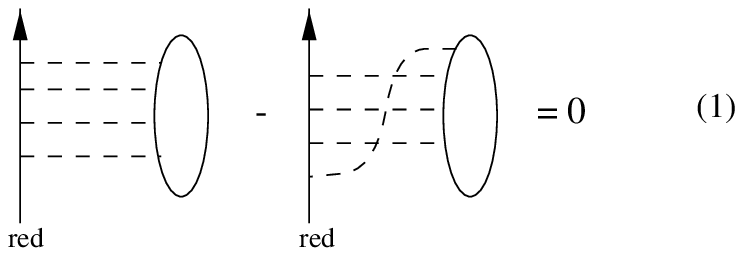}$$
    \caption{The link relation for chord diagrams} \label{F:link_rel}
    \end{figure}

\begin{figure}
    $$\psfig{file=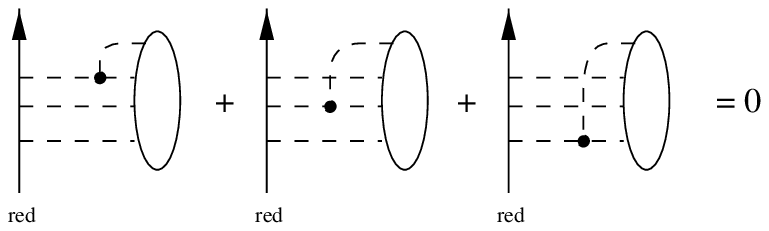}$$
    \caption{The link relation for Bounded Unitrivalent Diagrams} \label{F:link_rel2}
    \end{figure}

\begin{figure}
    $$\psfig{file=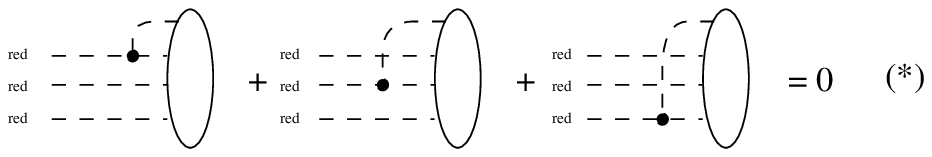}$$
    \caption{The link relation for unitrivalent diagrams} \label{F:link_rel3}
    \end{figure}

The obvious definition of chord diagrams for links is simply to replace the bounding line
segments with bounding circles.  However,
these diagrams are difficult to work with, and it is in particular unclear how to define the
unitrivalent diagrams.  Unlike for a knot, closing up the components of a string link of several
components is not a trivial operation, so we need to place some relations on the space
of unitrivalent diagrams.
Since we understand the spaces of chord diagrams and unitrivalent diagrams for string links, it
would be useful to be able to express these spaces for links as quotients of the spaces for
string links.  The question is then, what relations do we need?  One relation is fairly
obvious.  When we construct the space $A^l$ of bounded unitrivalent diagrams for links,
we replace the bounding line segments of $A^{sl}$ with directed circles.  Bar-Natan et. al.
observed (see Theorem 3, \cite{bgrt}) that this is exactly
equivalent to saying that the ``top'' edge incident to one of the line segments
can be brought around the circle to be on the ``bottom.''  So we can write $A^l$ as the
quotient of $A^{sl}$ by relation (1), shown in Figure~\ref{F:link_rel} (where the figure
shows {\it all} the chords with endpoints on the red component).
Then the Kontsevich integral for links, $Z^l$, is defined by cutting the link to make a string
link, applying the Kontsevich integral for string links, and then taking the quotient by
relation (1).  Now $w_d^l$ and $v_d^l$ are defined similarly to $w_d$ and $v_d$.
Given a link invariant $\gamma$ and a diagram $D$ in $A^l$, $w_d^l(\gamma)(D) =
\gamma(L_{\hat{D}})$, where $\hat{D}$ is the closure of the diagram $D$ (i.e. the bounding
line segments are closed to circles).  $L_{\hat{D}}$ is well-defined by Theorem 3 of
\cite{bgrt}.  Defining $v_d^l$ is even easier, now that we have $Z^l$.  Given a weight system
(element of the graded dual of $A^l$) $W$ and a link $L$, we define
$v_d^l(W)(L) = W(Z^l(L))$.
One advantage of this formulation of $A^l$ is that it enables us to define
the space $B^l$ of unitrivalent diagrams as a quotient of the already known
(see \cite{bn2}) space $B^{sl}$.  This was done by Bar-Natan et. al.~\cite{bgrt}.
Using the $STU$ relation, we can rewrite relation (1) as in Figure~\ref{F:link_rel2}.
This suggests how we should define the space $B^l$.  We will take the quotient of $B^{sl}$
by the relations (*) shown in Figure~\ref{F:link_rel3},
where the univalent vertices shown are {\it all} the univalent vertices of a given color.  With
these definitions, Bar-Natan et. al. proved that $A^l$ and $B^l$ are isomorphic:
\begin{thm} \label{T:isotopy-descends} (Theorem 3, \cite{bgrt})
The isomorphism between $A^{sl}$ and $B^{sl}$ descends to an isomorphism
between $A^l$ and $B^l$.
\end{thm}
\subsection{Link Concordance} \label{SS:concordance}
\begin{defn} \label{D:concordance}
Consider two k-component links $L_0$ and $L_1$, and two k-component string links $SL_0$
$SL_1$.  These can be thought of as embeddings (proper embeddings, in the case of the
string links):
$$L_i: \bigsqcup_{i=1}^k S^1 \hookrightarrow {\bf R}^3$$
$$SL_i: \bigsqcup_{i=1}^k I \hookrightarrow {\bf R}^2 \times I$$
A {\bf (link) concordance} between $L_0$ and $L_1$ is an embedding:
$$H: \left({\bigsqcup_{i=1}^k S^1}\right) \times I \hookrightarrow {\bf R}^3 \times I$$
such that $H(x,0) = (L_0(x),0)$ and $H(x,1) = (L_1(x),1)$.
Similarly, a {\bf (string link) concordance} between $SL_0$ and $SL_1$ is an embedding:
$$H: \left({\bigsqcup_{i=1}^k I}\right) \times I \hookrightarrow ({\bf R}^2 \times I) \times I$$
such that $H(x,0) = (SL_0(x),0),\ H(x,1) = (SL_1(x),1)$, and $H(i,t) = (i,t)$ for $i=0,1$.
A concordance is an isotopy if and only if H is level preserving; i.e. if the image of $H_t$ is
a (string) link at level $t$ for each $t \in I$.
\end{defn}
We want to extend the results of the last section to string links and links considered up to
concordance.  For string links, this has already been done by Habegger and Masbaum~\cite{hm}.
They describe the algebras $A^{csl}$ and $B^{csl}$ of bounded and unbounded unitrivalent diagrams
for string links up to concordance (which they denote $A^t$ and $B^t$), and observe that
they are isomorphic.  In brief, we take the quotient of $A^{sl}$
(resp. $B^{sl}$) by the space of diagrams with non-trivial
first homology.  In other words, we are left with tree diagrams.  It is then straightforward to
define $w_d^{csl}$ and $v_d^{csl}$ in the usual way, and show that they are
``almost'' inverses in the same sense that $w_d$ and $v_d$ are.
All of this extends to links just as it did for isotopy.  We define $A^{cl}$ as the quotient
of $A^{csl}$ by relation (1), and $B^{cl}$ as the quotient of $B^{csl}$ by relation (*).  We
then define $Z^{cl}$, $w_d^{cl}$, and $v_d^{cl}$ just as we did for links up to homotopy.
Finally, the arguments of Bar-Natan et. al. carry through to show:
\begin{thm} \label{T:concordance-descends} (Theorem 3, \cite{bgrt})
The isomorphism between $A^{csl}$ and $B^{csl}$ descends to an isomorphism
between $A^{cl}$ and $B^{cl}$.
\end{thm}
{\bf Remark:}  By results of Habegger and Masbaum (see Theorem 5.5 of \cite{hm}), $Z^{cl}$ is
the {\it universal} finite type invariant of link concordance.  By this we mean that it dominates
all other such invariants.
\section{The Size of $B^{cl}$} \label{S:size}
Now that we have properly defined the space $B^{cl}$ of unitrivalent diagrams for link homotopy, we want
to analyze it more closely.  We will consider the case when $B^{cl}$ is a vector space over the reals
(or, more generally, a module over a ring of characteristic 0).
In particular, we would like to know exactly which diagrams of $B^{csl}$
are in the kernel of the relation (*) (i.e. are 0 modulo (*)).  We will find that the answer is
``almost everything'' - to be precise, any unitrivalent diagram with a component of degree 2 or more.
We will start by proving a couple of base cases, and then prove the rest of the theorem by induction.
Let $B^{csl}(k)$ denote the space of unitrivalent diagrams for string link concordance
with $k$ possible colors for the univalent vertices (i.e. we are looking at links with $k$ components).
Consider a diagram $D \in B^{csl}(k)$.  Recall from the previous sections that each component of $D$
is a tree diagram.
{\bf Notation:}  Before we continue, we will introduce two bits of notation which will be useful in this
section.
\begin{itemize}
    \item  Given a unitrivalent diagrams $D$, we define $m(D;i,j)$ to be the number of
components of $D$ which are simply line segments with ends colored $i$ and $j$, as shown below:
$$i-----j$$
    \item  Components of a diagram with degree greater than one will be called {\it large}
components.  Components of degree one will be called {\it small} components.
\end{itemize}
\subsection{Knots and Two Component Links} \label{SS:1&2comp}
\begin{figure}
    $$\psfig{file=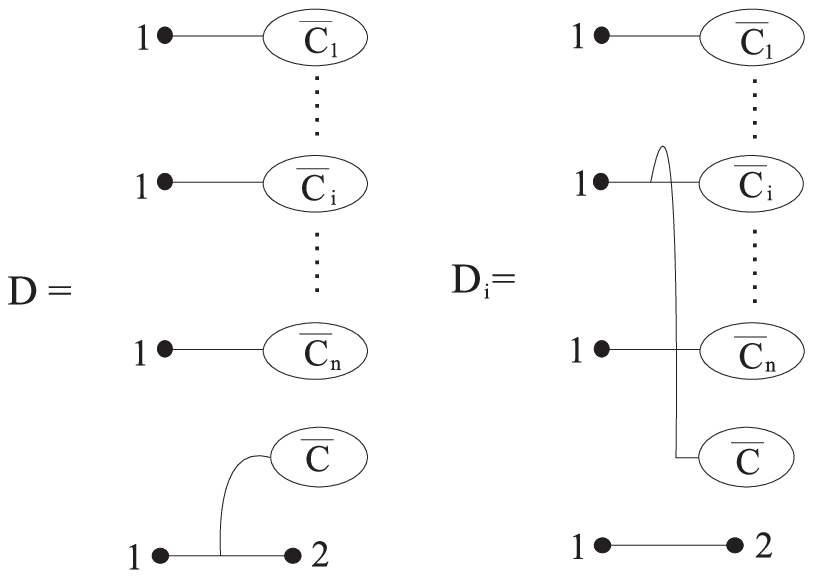}$$
    \caption{Diagrams arising from relation (*)} \label{F:arise}
    \end{figure}

We begin by considering the case of knot concordance, when $k = 1$.  Ng~\cite{ng} has already shown
that the only finite type invariant of knot concordance is the ${\bf Z}_2$-valued Arf invariant, so
there are no real-valued finite type invariants of knot concordance.  We will begin, as a
warm-up, by showing this result using unitrivalent diagrams.  Bar-Natan has shown that the spaces of
unitrivalent diagrams for knots and string links of one component are isomorphic, so the relation
(*) has no effect.
\begin{lem} \label{L:knots}
$B^{csl}(1) = B^{cl}(1) = 0.$
\end{lem}
{\sc Proof:}  First, we consider the case when $D \in B^{csl}(1)$ has a large component $C$.
Since $C$ is a tree, we can use the $IHX$ relation as in Figure~\ref{F:treepf}
to rewrite $C$ as a sum of diagrams:
    \begin{figure}
    $$\psfig{file=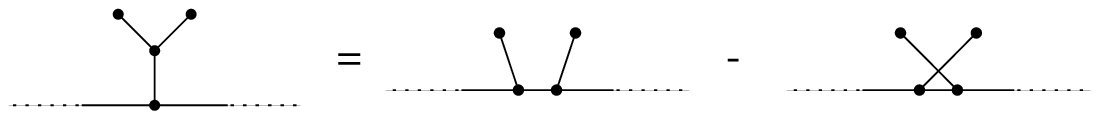}$$
    \caption{Using IHX to get ``standard'' tree diagram} \label{F:treepf}
    \end{figure}
$$\psfig{file=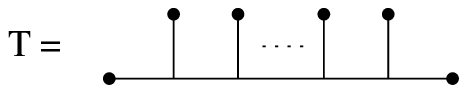}$$
Clearly, twisting the branch at the right or left end of $T$ gives us the same diagram, but by
the $AS$ relation, this flips the sign.  Therefore, $T = -T = 0$, so $D = 0$.
The case when all the components of D have degree one is somewhat more subtle, though it
is essentially an application of the 1-term relation.  If D has only small components,
then $\chi(D)$ (see Section~\ref{SS:ud}) is a linear combination of bounded unitrivalent diagrams
with no interior vertices (i.e. chord diagrams).  By repeated applications of the $STU$ relation,
we can isolate a chord in each of these diagrams, at the expense of adding a linear combination
of diagrams which {\it do} have internal vertices.  The diagrams with isolated chords disappear
by the 1-term relation, so we are left with a linear combination of bounded unitrivalent diagrams
with at least one internal vertex.  Now we apply $\sigma$ to this linear combination to get a
linear combination of (unbounded) unitrivalent diagrams.  Since the basic operations $U$ and $S$
of $\sigma$ can never decrease the number of internal vertices, every unitrivalent diagram in
the image of $\sigma$ will have at least one internal vertex; i.e. at least one large component.
Hence, by the first case, all of these diagrams are 0.  So we have
shown that $\sigma(\chi(D)) = 0$.  But $\sigma$ and $\chi$ are inverse isomorphisms, so this means
$D = 0$, as desired.  $\Box$
Next we consider links of two components, i.e. $k=2$.  In this case we need to consider the
effect of the relation (*).
\begin{lem} \label{L:k=2}
Let $D \in B^{csl}(2)$.  If $D$ has a component $C$ of degree $d \geq 2$, then $D$ is trivial
modulo (*).
\end{lem}
{\sc Proof:}  $D$ is a diagram with all endpoints colored 1 or 2.  Note that $C$ must have endpoints
of both colors, or $D$ will be trivial by the same argument as in Lemma~\ref{L:knots}.  In fact,
any terminal branch of $C$ must have the form (where $\bar{C}$ denotes the remainder of $C$):
$$C:\ \ \begin{matrix} \bar{C} \\ | \\ | \\ 1-----2 \end{matrix}$$
Otherwise, if the two endpoints have the same color, $C$ (and hence $D$) is trivial by
the $AS$ relation.
The proof is by induction on the number of large components of $D$.  In
the base case, there is only one such component, $C$.  So all the other components $C_i$ are
simply line segments labeled $a$ and $b$, where $a,b \in \{1,2\}$, as shown:
$$C_i:\ \ a-----b$$
Now we apply the relation (*) to the branch of $C$ shown above using the color 1, as in
Figure~\ref{F:arise} (where $\bar{C_i}$ represents the remainder of the component $C_i$).
We denote the merger of $C_i$ and $\bar{C}$ as $C'_i$.  We do not need to consider other
vertices of $C$ colored 1, since these terms result in a diagram with a loop, which are trivial
in concordance.  Since we are expanding using the color 1, we only need to consider $C_i$ where
$a$ or $b$ is 1.  If they are both 1, then by the $AS$ relation:
$$C'_i:\ \ \begin{matrix} \bar{C} \\ | \\ | \\ 1-----1 \end{matrix} = 0$$
Therefore, we need only consider $C_i$ with one endpoint colored 1 and the other colored 2.
But in this case, $C'_i = C$, so $D_i = D$, and we get an equation $D + \sum{D_i} = (1+m)D = 0$
for some $m \geq 0$.  Hence $D = 0$ modulo (*).  This concludes the base case.
For the inductive step, we assume that the theorem is true for diagrams with $n$ large
components, and consider a diagram $D$ with $n+1$ large components.  Let $C$ be one
of these components.  As before, C has a branch:
$$C:\ \ \begin{matrix} \bar{C} \\ | \\ | \\ 1-----2 \end{matrix}$$
Also, every small component $C_i$ looks like (where $a,b \in \{1,2\}$):
$$C_i:\ \ a-----b$$
Once again, we apply the relation (*) to the color 1, using this branch of $C$, and find that
$D + \sum{D_i} = 0$.  Whenever $C_i$ is small, $D_i = D$ (as in the base case, we need only
consider the case when $a = 1$ and $b = 2$).  If $C_i$ is large, then $D_i$ has
one fewer large component than $D$ does,
since the bulk of $C$ has joined with $C_i$, leaving behind a line segment of degree one.  So
$D_i$ is trivial modulo (*) by the inductive hypothesis.  Hence we are left with a sum of
copies of $D$, and conclude that $(1+m)D = 0$ for some $m \geq 0$, so $D = 0$ modulo (*).
This concludes the induction and the proof.  $\Box$
\subsection{Three Component Links} \label{SS:3comp}
\begin{figure}
    $$\psfig{file=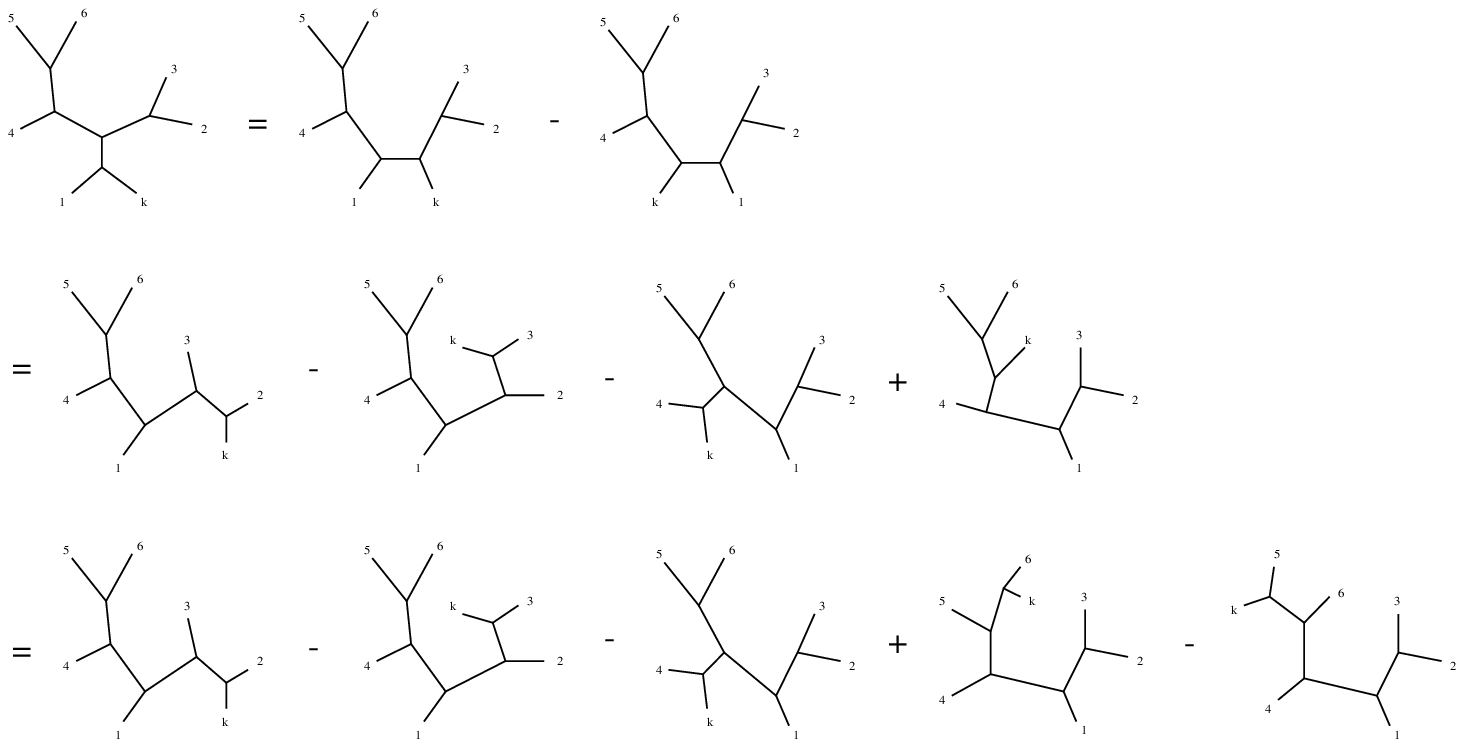}$$
    \caption{Using the IHX relation to decompose a diagram} \label{F:expand}
    \end{figure}

The case when $k=3$ is our final special case before the proof of the general theorem, and
it is significantly more complicated than the previous two lemmas.  The main step is to
show that no component of $D$ can have two endpoints of the same color (essentially, this
reduces the problem to the case of link homotopy, treated in \cite{me1}).  Once again, we will be
inducting on the number of large components of the diagram.  The base
case contains most of the work of the proof, so we present it as a separate lemma.
\begin{lem} \label{L:k=3_base}
If k=3, D has exactly one large component C, and C has two endpoints of the
same color, then D is trivial modulo (*).
\end{lem}
{\sc Proof:}  Without loss of generality, we will say that $C$ has two endpoints colored 1.
If these endpoints are on the same final branch, as shown below, then $D$ will be trivial by
the $AS$ relation (since we will have $D = -D$).  $\bar{C}$ denotes the remainder of C:
$$C:\ \ \begin{matrix} \bar{C} \\ | \\ | \\ 1-----1 \end{matrix}$$
Otherwise, we can use the $IHX$ relation to move one of the endpoints colored 1 out to the
ends of the component $C$, as shown in Figure~\ref{F:expand} (where we move the
endpoint colored $k$).
So it suffices, without loss
of generality, to consider $C$ with a branch as shown below:
$$C:\ \ \begin{matrix} \bar{C} \\ | \\ | \\ 1-----2 \end{matrix}$$
Our proof will be by induction on $m(D;2,3)$, inducting among diagrams with exactly one
large component, which have a branch with colors 1 and 2 as shown above.
For the base case, $m(D;2,3) = 0$.  This means that there are no small components of $D$
colored 2 and 3.  So the only degree 1 components with an endpoint colored 2
have their other endpoint colored 1 or 2:
$$C_i:\ \ a-----2,\ a\in \{1,2\}$$
Now we apply the relation (*) to $D$ as in Figure~\ref{F:arise}, only now we are applying it
using the color 2.  So we only need to consider the components of $D$ with endpoints colored
2.  Because of the loop relation for concordance, we can ignore other endpoints of $C$ which
might be colored 2, and just consider the remaining components.  These are all small
components with endpoints colored $a$ and 2 as described above.  As in Lemma~\ref{L:k=2},
the case when $a=2$ can be ignored by the $AS$ relation, so we are reduced to the case of
small components with endpoints colored 1 and 2.  Therefore, we obtain the
equation $D + \sum{D_i} = 0$, where every $D_i = D$.  So $D + m(D;1,2)D =
(1+m(D;1,2))D = 0$.  Since $m(D;1,2) \geq 0$, we conclude that $D = 0$ modulo (*).  This
concludes the base case of the induction.
For the inductive step, we will have small components with endpoints colored 2 and 3.
Our goal is to reduce the number of such components.  Now, when we apply the relation (*),
we can again ignore small components with both endpoints colored 2, by the $AS$ relation.
We get an equation $D + m(D;1,2)D + m(D;2,3)D' = 0$, where $D'$ is identical to $D$ except
that a small component with endpoints colored 2 and 3 has been replaced by one with
endpoints colored 1 and 2, and the endpoint of $C$ colored 1 above has been colored 3.  We
denote this analogue of $C$ in $D'$ by $C'$.  Diagramatically, we can represent $D'$ by
showing the changes that have been made to $D$:
$$D':\ \ \begin{matrix} \bar{C} \\ | \\ | \\ 3-----2 \end{matrix} \  (2,3) \rightarrow (1,2)$$
Note that $m(D';2,3) = m(D;2,3) - 1$.
As we did in Lemma~\ref{L:knots}, we can apply the $IHX$ relation to $C'$, fixing the branch
shown above.  As a result, it suffices to consider the case when $C'$ has the form shown
below (where $n$ is the degree of $C'$):
$$\psfig{file=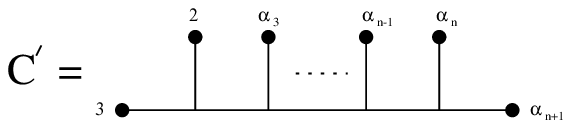}$$
We can assume that one of the endpoints $\alpha_3,...,\alpha_n$ is colored 1.  If not, then
$\alpha_{n+1} = 1$ (since we know there is a second endpoint of $C$ colored 1), and we can
switch $\alpha_n$ and $\alpha_{n+1}$ using the $AS$ relation (at the cost of reversing the
sign of $D'$).
We can apply the $IHX$ relation as in Figure~\ref{F:treepf} to move the endpoint
colored 2 along the ``spine'' of $C'$
(i.e. the path from the endpoint colored 3 to the endpoint colored $\alpha_{n+1}$) and
obtain the decomposition $D' = \sum{D^i}$, where we transform $C$ into $C^i$ as shown below:
$$\psfig{file=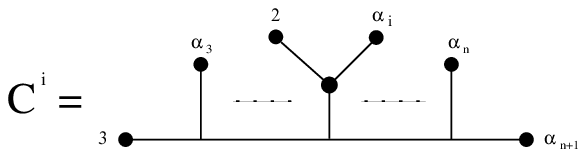}$$
Notice that the spine of $C^i$ is shorter than the spine of $C$ by one edge; and that one
of the branches along the spine has grown correspondingly.  We observe that if $\alpha_i = 1$,
then $D^i = 0$ by the inductive hypothesis, since it now has a branch colored 1 and 2, and
$m(D^i;2,3) = m(D';2,3) = m(D;2,3) - 1$.  So the endpoint colored 1 will not be incorporated
into the larger branch, and has moved one position closer to the endpoint colored 3 at the
far left.
In the case when $\alpha_3 = 1$, $D^3 = 0$, and $C^i$ has the form below for $i > 3$:
$$\psfig{file=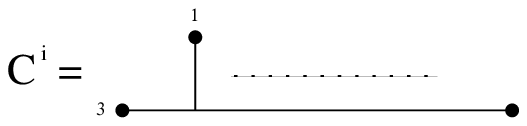}$$
Now if we use the relation (*) to expand $D^i$ along the color 1, we find that
$D^i + m(D^i;1,3)D^i + m(D^i;1,2)D^i_2 = 0$, where $D^i_2$ is the result of replacing a
small component colored 1 and 2 with a small component colored 1 and 3, and changing
the endpoint of $C^i$ colored 3 to one colored 2.  Then $D^i_2$ has a large component
with a branch colored 1 and 2, and $m(D^i_2;2,3) = m(D^i;2,3) = m(D;2,3)-1$,
so $D^i_2$ is trivial modulo (*) by the inductive hypothesis.  Therefore, we find that
$(1+m(D^i;1,3))D^i = 0$, and hence $D^i = 0$ modulo (*).
If $\alpha_3 \neq 1$, the first branch on the spine of $C^i$ (adjacent to the endpoint
colored 3) will only have endpoints colored 2 and 3 (not 1).  Now we repeat the process
by expanding the (non-trivial) $C^i$'s, using the first branch on the spine of $C^i$.
In order to continue the process, we need to show the following fact:
{\sc Claim:}  If $K$ is a unitrivalent diagram with all endpoints colored 2 or 3, then
the following diagram is trivial modulo (*):
$$\psfig{file=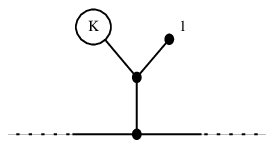}$$
Here $K$ is assumed to be a subdiagram of the only component of degree greater than 1
in a diagram $E$ such that $m(E;2,3) = m(D;2,3)-1$.
{\sc Proof of Claim:}  Using the $IHX$ relation as in Figure~\ref{F:expand}, we can
decompose this diagram into a linear combination of diagrams where the endpoint
colored 1 has migrated out to one of the ends of $K$, leaving a diagram with a branch
with endpoints colored 1 and 2 or 1 and 3.  In the first case, the diagram is trivial
by the inductive hypothesis.  In the second case, the diagram is trivial by the argument
used above for $\alpha_3 = 1$.  So we conclude that $E$ is trivial modulo (*).
$\Box$ (Claim)
Using the claim, we can continue the process, moving the endpoint colored 1 to the left
at each stage, until we are left with diagrams where the
first branch, adjacent to the endpoint colored 3, consists of a single endpoint colored
1.  These diagrams are trivial by the argument above (when $\alpha_3 = 1$).  We can
conclude that all of the $D^i$'s are trivial modulo (*).  Therefore, $D'$ is trivial
modulo (*), and we obtain the equation $D + m(D;1,2)D = (1+m(D;1,2))D = 0$.  We conclude
that $D = 0$ modulo (*), which finishes the induction and the proof. $\Box$
\begin{lem} \label{L:k=3_induct}
If k=3 and D has a component C with two endpoints of the same color, then D is trivial
modulo (*).
\end{lem}
{\sc Proof:}  First, we assume that $C$ is large.  We will induct
on the number of large components of $D$.  The base case of the
induction has already been proved, in Lemma~\ref{L:k=3_base}.  The general case follows
exactly the same argument.  The only modification is to notice that, whenever (*) is
applied, the diagrams $D_i$ (as in Figure~\ref{F:arise}) which arise from large components
$C_i$ will have fewer large components than
D (since $C$ and $C_i$ have been joined, leaving behind a component of degree one).  So
by the inductive hypothesis, these diagrams can be ignored at every stage.  Therefore,
exactly the same proof shows that $D$ is trivial modulo (*).
If $C$ is small, with both endpoints the same color, we can use the previous case
together with the argument from Lemma~\ref{L:knots} to show $D$ is trivial.  $\Box$
\begin{lem} \label{L:k=3}
If k=3 and D has a large component C, then D is trivial modulo (*).
\end{lem}
{\sc Proof:}  By Lemma~\ref{L:k=3_induct}, there is (up to sign) only one possible
diagram for $C$:
$$C = \ \begin{matrix} 3 \\ | \\ | \\ 1-----2 \end{matrix}$$
Now we apply the relation (*) to $D$ using C and the color 1.  By Lemma~\ref{L:k=3_induct},
we need only consider components with endpoints colored 1 or 2, and no component can have
two endpoints of the same color.  Therefore, we need only consider components $C_i$ as
shown:
$$C_i = \ 1-----2$$
This gives us the equation $D + m(D;1,2)D = (1+m(D;1,2))D = 0$.  We conclude that $D = 0$
modulo (*), which completes the proof for the case $k = 3$.  $\Box$
\subsection{The General Case} \label{SS:general}
We are now ready to begin our proof of the general case.  As we did for the case when
$k = 3$, we will induct on the number of large components of $D$.  Once again,
for clarity, we will prove the base case (which contains most of the work) as a
separate lemma.
\begin{lem} \label{L:base}
If D has exactly one large component C, then D is trivial modulo (*).
\end{lem}
{\sc Proof:}  The method of proof for this lemma is very similar to the earlier lemmas.
We will successively apply (*) (and do a single expansion via IHX)
until we obtain a set of diagrams which are all either trivial or repetitions of earlier
diagrams.  We can then backtrack to show that everything disappears.  However, we will
need to apply (*) four times.  This
unfortunately makes keeping track of the diagrams somewhat confusing - we have done our best.
Without loss of generality, as before, we can assume that $C$ has a branch as shown:
$$C:\ \ \begin{matrix} \bar{C} \\ | \\ | \\ 1-----2 \end{matrix}$$
We apply (*) using the color 1 and find that $D + m(D;1,2)D + \sum_{a\neq 1,2}{m(D;1,a)D_a} = 0$,
where $D_a$ is the same as $D$ except that:
\begin{itemize}
    \item  $C$ has been replaced by a component $C_a$ identical to it except that the
endpoint colored 2 is now colored $a$ (so $\bar{C_a} = \bar{C}$)
    \item  A line segment with endpoints colored 1 and $a$ has been replaced by a line segment
with endpoints colored 1 and 2.  In other words, $m(D_a;1,a) = m(D;1,a)-1$ and $m(D_a;1,2) =
m(D;1,2)+1$.
\end{itemize}
We will denote this as shown below:
$$D_a:\ \ \begin{matrix} \bar{C} \\ | \\ | \\ 1-----a \end{matrix} \  (1,a) \rightarrow (1,2)$$
As in Lemma~\ref{L:k=3_base}, we use the IHX relation to decompose $D_a =
\sum_{\alpha_i \neq a}{\pm D_a^i}$ (summing over the endpoints of $C_a$, with colors $\alpha_i$),
where the analogue $C_a^i$ of $C_a$ in $D_a^i$ has a branch as shown, and the other components of the
diagram are the same as $D_a$:
$$D_a^i:\ \ \begin{matrix} \bar{C_a^i} \\ | \\ | \\ i-----a \end{matrix} \ (1,a) \rightarrow (1,2)$$
By abuse of notation, we write the color $\alpha_i$ as simply $i$.  Since we will never have to
compare different $D_a^ii$'s, this will not cause any confusion.
Note that, aside from having endpoints of the same colors, $C_a^i$ looks nothing like $C_a$.
Now we apply (*) to $D_a^i$, using the color $i$.  In the pictures we use to describe the various
diagrams that we produce in what follows, we will just be showing how the diagrams differ from
$D_a^i$.  This will involve showing how $C_a^i$ has been altered, and which line segments have
been added or removed.  At each stage, we will eliminate loop diagrams without comment.  We
obtain the relation:
$$D_a^i + m(D_a^i;i,a)D_a^i + m(D_a^i;2,i)D_{a2}^i +
\sum_{b\neq i,2,a}{m(D_a^i;i,b)D_{ab}^i} = 0$$
where:
$$D_{a2}^i:\ \ \begin{matrix} \bar{C_a^i} \\ | \\ | \\ i-----2 \end{matrix} \ (2,i) \rightarrow (i,a)$$
$$D_{ab}^i:\ \ \begin{matrix} \bar{C_a^i} \\ | \\ | \\ i-----b \end{matrix} \ (i,b) \rightarrow (i,a)$$
Next we apply (*) to $D_{ab}^i$, using the color $b$, and find that:
$$D_{ab}^i + m(D_{ab}^i;i,b)D_{ab}^i + m(D_{ab}^i;2,b)D_{ab2}^i +
\sum_{c\neq b,i,2}{m(D_{ab}^i;b,c)D_{abc}^i} = 0$$
where:
$$D_{ab2}^i:\ \ \begin{matrix} \bar{C_a^i} \\ | \\ | \\ 2-----b \end{matrix} \
\begin{matrix} (i,b) \rightarrow (i,a) \\ (2,b) \rightarrow (i,b) \end{matrix}
\Rightarrow (2,b) \rightarrow (i,a)$$
$$D_{abc}^i:\ \ \begin{matrix} \bar{C_a^i} \\ | \\ | \\ c-----b \end{matrix} \
\begin{matrix} (i,b) \rightarrow (i,a) \\ (c,b) \rightarrow (i,b) \end{matrix}
\Rightarrow (c,b) \rightarrow (i,a)$$
Now we apply (*) to $D_{ab2}^i$ using the color 2, and to $D_{abc}^i$, using the color $c$.  We
get two relations:
$$D_{ab2}^i + m(D_{ab2}^i;2,b)D_{ab2}^i + m(D_{ab2}^i;2,i)D_{ab2i}^i +
\sum_{c\neq b,i,2}{m(D_{ab2}^i;2,c)D_{ab2c}^i} = 0$$
$$D_{abc}^i + m(D_{abc}^i;b,c)D_{abc}^i + m(D_{abc}^i;2,c)D_{abc2}^i + m(D_{abc}^i;i,c)D_{abci}^i +$$
$$\sum_{d\neq b,c,i,2}{m(D_{abc}^i;c,d)D_{abcd}^i} = 0$$
where:
$$D_{ab2i}^i:\ \ \begin{matrix} \bar{C_a^i} \\ | \\ | \\ 2-----i \end{matrix} \
\begin{matrix} (2,b) \rightarrow (i,a) \\ (2,i) \rightarrow (2,b) \end{matrix}
\Rightarrow (2,i) \rightarrow (i,a)$$
$$D_{ab2c}^i:\ \ \begin{matrix} \bar{C_a^i} \\ | \\ | \\ 2-----c \end{matrix} \
\begin{matrix} (2,b) \rightarrow (i,a) \\ (2,c) \rightarrow (2,b) \end{matrix}
\Rightarrow (2,c) \rightarrow (i,a)$$
$$D_{abc2}^i:\ \ \begin{matrix} \bar{C_a^i} \\ | \\ | \\ c-----2 \end{matrix} \
\begin{matrix} (c,b) \rightarrow (i,a) \\ (2,c) \rightarrow (c,b) \end{matrix}
\Rightarrow (2,c) \rightarrow (i,a)$$
$$D_{abci}^i:\ \ \begin{matrix} \bar{C_a^i} \\ | \\ | \\ c-----i \end{matrix} \
\begin{matrix} (c,b) \rightarrow (i,a) \\ (i,c) \rightarrow (c,b) \end{matrix}
\Rightarrow (i,c) \rightarrow (i,a)$$
$$D_{abcd}^i:\ \ \begin{matrix} \bar{C_a^i} \\ | \\ | \\ c-----d \end{matrix} \
\begin{matrix} (c,b) \rightarrow (i,a) \\ (c,d) \rightarrow (c,b) \end{matrix}
\Rightarrow (c,d) \rightarrow (i,a)$$
We make several observations:
\begin{itemize}
    \item  $D_{ab2i}^i = -D_{a2}^i$.
    \item  $D_{ab2c}^i = -D_{abc2}^i$.
    \item  $D_{abci}^i = -D_{ac}^i$.
    \item  $D_{abcd}^i = D_{adc}^i = -D_{acd}^i$.
\end{itemize}
Now that we have these recursive relations, we can plug them into our various equations.
We will use the following equalities:
\begin{eqnarray*}
m(D_{ab}^i;i,b) + 1 & = & m(D_a^i;i,b) \\
m(D_{ab2}^i;2,b) + 1 & = & m(D_{ab}^i;2,b) = m(D_a^i;2,b) \\
m(D_{abc}^i;b,c) + 1 & = & m(D_a^i;b,c)
\end{eqnarray*}
And for all the other coefficients we have:
$$m(D_*^i;x,y) = m(D_a^i;x,y)$$
For convenience, we will write $m(x,y) = m(D_a^i;x,y)$ in what follows:
\begin{eqnarray*}
(m(i,a)+1)D_a^i + m(2,i)D_{a2}^i & = & \sum_{b\neq a,i,2}{-m(i,b)D_{ab}^i} \\
 & = & \sum_{b\neq a,i,2}{-(m(D_{ab}^i;i,b)+1)D_{ab}^i} \\
 & = & \sum_{b\neq a,i,2}{\left({m(2,b)D_{ab2}^i +
    \sum_{c\neq b,i,2}{m(b,c)D_{abc}^i}}\right)}
\end{eqnarray*}
Note that:
\begin{eqnarray*}
m(2,b)D_{ab2}^i & = & (m(D_{ab2}^i;2,b)+1)D_{ab2}^i \\
 & = & -m(2,i)D_{ab2i}^i + \sum_{c\neq b,i,2}{-m(2,c)D_{ab2c}^i} \\
m(b,c)D_{abc}^i & = & (m(D_{abc}^i;b,c)+1)D_{abc}^i \\
 & = & -m(2,c)D_{abc2}^i - m(i,c)D_{abci}^i - \sum_{d\neq c,b,i,2}{m(c,d)D_{abcd}^i}
\end{eqnarray*}
Therefore:
$$m(2,b)D_{ab2}^i + \sum_{c\neq b,i,2}{m(b,c)D_{abc}^i} = $$
$$-m(2,i)D_{ab2i}^i + \sum_{c\neq b,i,2}{\left({-m(2,c)(D_{ab2c}^i + D_{abc2}^i)
    - m(i,c)D_{abci}^i - \sum_{d\neq c,b,i,2}{m(c,d)D_{abcd}^i}}\right)} = $$
$$m(2,i)D_{a2}^i + \sum_{c\neq b,i,2}{\left({m(i,c)D_{ac}^i +
    \sum_{d\neq c,b,i,2}{m(c,d)D_{acd}^i}}\right)}$$
We plug this back in above to find:
$$(m(i,a)+1)D_a^i + m(2,i)D_{a2}^i = \sum_{b\neq a,i,2}{\left({m(2,i)D_{a2}^i +
    \sum_{c\neq b,i,2}{\left({m(i,c)D_{ac}^i +
    \sum_{d\neq c,b,i,2}{m(c,d)D_{acd}^i}}\right)}}\right)}$$
We notice that:
\begin{eqnarray*}
\sum_{c\neq b,i,2\ }{\sum_{d\neq c,b,i,2}{m(c,d)D_{acd}^i}} & = &
    \frac{1}{2}\sum_{\substack{c\neq d \\ c,d\neq b,i,2}}{m(c,d)(D_{acd}^i + D_{adc}^i)} \\
 & = & \frac{1}{2}\sum_{\substack{c\neq d \\ c,d\neq b,i,2}}{m(c,d)(D_{acd}^i - D_{acd}^i)} \\
 & = & 0
\end{eqnarray*}
Therefore:
$$(m(i,a)+1)D_a^i + m(2,i)D_{a2}^i = \sum_{b\neq a,i,2}{\left({m(2,i)D_{a2}^i +
    \sum_{c\neq b,i,2}{m(i,c)D_{ac}^i}}\right)}$$
Returning to our first equation, we obtain (simply replacing $b$ by $c$ in the second
equality):
\begin{eqnarray*}
(m(i,a)+1)D_a^i + m(2,i)D_{a2}^i & = & \sum_{b\neq a,i,2}{-m(i,b)D_{ab}^i} \\
 & = & \sum_{c\neq a,i,2}{-m(i,c)D_{ac}^i} \\
 & = & \left({\sum_{c\neq b,i,2}{-m(i,c)D_{ac}^i}}\right) + m(i,a)D_{aa}^i - m(i,b)D_{ab}^i
\end{eqnarray*}
Since $D_{aa}^i = D_a^i$, we can cancel and rearrange terms to write:
$$\sum_{c\neq b,i,2}{m(i,c)D_{ac}^i} = -D_a^i - m(i,b)D_{ab}^i - m(2,i)D_{a2}^i$$
Therefore:
\begin{eqnarray*}
(m(i,a)+1)D_a^i + m(2,i)D_{a2}^i & = & \sum_{b\neq a,i,2}{\left({m(2,i)D_{a2}^i - D_a^i -
    m(i,b)D_{ab}^i - m(2,i)D_{a2}^i}\right)} \\
 & = & \sum_{b\neq a,i,2}{\left({-D_a^i - m(i,b)D_{ab}^i}\right)} \\
 & = & -\left({\sum_{b\neq a,i,2}{D_a^i}}\right) + (m(i,a)+1)D_a^i + m(2,i)D_{a2}^i
\end{eqnarray*}
Finally, cancelling gives us:
$$\sum_{b\neq a,i,2}{D_a^i} = 0$$
We now need to consider the possibilities for the color $i$.  If $i = 2$, then:
$$\sum_{b\neq a,i,2}{D_a^i} = \sum_{b\neq a,2}{D_a^i} = (k-2)D_a^i = 0$$
However, if $i \neq 2$, then:
$$\sum_{b\neq a,i,2}{D_a^i} = (k-3)D_a^i = 0$$
Since we have already dealt with the cases when $k \leq 3$ in Lemma~\ref{L:knots},
Lemma~\ref{L:k=2} and Lemma~\ref{L:k=3}, we can conclude that
$D_a^i$ is trivial modulo (*) for every $i$.  Hence, $D_a$ is trivial for every $a$.  Finally,
this means that $(1+m(D;1,2))D = 0$, so $D$ will also be trivial modulo (*).  This completes the
proof.  $\Box$
\begin{thm} \label{T:induction}
If D has any large components C, then D is trivial modulo (*).
\end{thm}
{\sc Proof:}  As in Lemma~\ref{L:k=3_induct}, we induct on the number of large components
of $D$, and the argument for the general case is almost identical to the
argument for the base case.  We just have to observe that whenever we apply the relation (*),
we can ignore large components of $D$ (other than the one we're working with),
since the diagrams they give rise to have fewer large components than
$D$ does.  This completes the induction and the proof.  $\Box$
This theorem tells us that the only elements of $B^{csl}$ which are {\it not} in the kernel of the
relation (*) are unitrivalent diagrams all of whose components are of degree 1 (i.e. line segments).
By the arguments of Lemma~\ref{L:knots}, we also know that the degree 1 components with both
endpoints of the same color are also trivial, so we only need to consider line segments with
different colors on the two endpoints.  Restricted to the space generated by these elements,
(*) is clearly trivial, so $B^{cl}$ is in fact simply the polynomial algebra over the
reals generated by these unitrivalent diagrams (since (*) is trivial on this space,
$B^{cl}$ inherits a multiplication from $B^{csl}$).  We formalize this as a corollary:
\begin{cor} \label{C:algebra}
$B^{cl}(k)$ (and hence $A^{cl}(k)$) is isomorphic to the algebra
${\bf R}[x_{ij}]$,
where each $x_{ij}$ is of degree 1, and $1\leq i < j\leq k$.
\end{cor}
It is well-known that these diagrams correspond to the pairwise
linking numbers of the components, so we conclude:
\begin{cor} \label{C:reduction}
The pairwise linking numbers of the components of a link are the only finite type link concordance
invariants of the link.
\end{cor}

\section{Acknowledgements}
I wish to thank my advisor, Robion Kirby, for his advice and support, and all the
members of the Informal Topology Seminar at Berkeley.

\end{document}